\documentclass[reqno,fleqn]{amsart}%
\usepackage{amsmath}
\usepackage{amssymb}
\usepackage{amsfonts}
\usepackage{graphicx}
\usepackage{amsthm}
\usepackage{enumerate}

\newcommand{\R}{\mathbb{R}}
\newcommand{\C}{\mathbb{C}}
\newcommand{\f}{\rightarrow}
\newcommand{\deb}{\bar\partial}
\newcommand{\de}{\partial}

\newcommand{\K}{K\"{a}hler}
\newcommand{\ngh}{neighbourhood}
\newcommand{\M}{\mathcal M}
\newcommand{\lmb}{\lambda}

\newcommand{\E}{\operatorname{Exp}}
\newcommand{\hyp}{\operatorname{hyp}}

\newcommand{\tr}{\operatorname{tr}}

\newcommand{\id}{\operatorname{id}}
\newcommand{\Cut}{\operatorname{Cut}}

\newtheorem{thm}{Theorem}
\newtheorem{prop}{Proposition}

\newtheorem{exmp}[prop]{Example}

\begin{document}
\title{The diastatic exponential of  a symmetric space}
\author[A. Loi]{Andrea Loi}
\address{Dipartimento di Matematica e Informatica, Universit\`{a} di Cagliari,
Via Ospedale 72, 09124 Cagliari, Italy}
\email{loi@unica.it;   roberto.mossa@gmail.com }
\author[R. Mossa]{Roberto Mossa}
\thanks{Research partially supported by GNSAGA (INdAM) and MIUR of Italy}
\date{April 7, 2009}

\begin{abstract}
Let $\left(M, g\right)$ be a real analytic \K\ manifold.
We say that a smooth map
$\E_p:W\rightarrow M$
from a \ngh\  $W$ of the origin of $T_pM$ into $M$
is a {\em diastatic exponential} at $p$
if it satisfies
$$
\left(d \E_p\right)_0=\id_{T_pM},
$$
$$
D_p\left( \E_p \left(v\right)\right)=g_p\left(v, v\right),\  \forall v\in W,
$$
where $D_p$ is Calabi's diastasis function at $p$
(the usual exponential $\exp_p$ obviously satisfied these
equations when $D_p$ is replaced by the square of the geodesics distance $d^2_p$
from $p$).
In this paper we prove that for every point $p$ of an Hermitian symmetric space of noncompact type M there exists a globally defined diastatic exponential centered in $p$ which is a diffeomorphism
and it is uniquely determined by its restriction to polydisks.
An analogous result   holds true
in an open  dense  \ngh\  of every point of $M^*$, the compact dual of $M$ .
We also provide a geometric  interpretation
of the symplectic duality map  (recently introduced  in \cite{sympdual})  in terms of diastatic exponentials.
As a byproduct  of our analysis we show that the symplectic duality  map
pulls back the reproducing  kernel of $M^*$
to the reproducing kernel of $M$.
\end{abstract}
\maketitle

\section*{Introduction and statements of the main results}
Let $M$ be a  $n$-dimensional complex  manifold endowed with a real analytic    \K\ metric $g$.
For a fixed point $p\in M$ let  $D_p:U\rightarrow \R$  be the  {\em Calabi diastasis function},
defined  in the following way. Recall that
a \K\ potential is an
analytic function  $\Phi$
defined in a neighborhood
of a point $p$ such that
$\omega =\frac{i}{2}\partial\bar\partial\Phi$,
where $\omega$ is the \K\ form associated
to $g$.
By duplicating the variables $z$ and $\bar z$
a potential $\Phi$ can be complex analytically
continued to a function
$\tilde\Phi$ defined in a neighborhood
$U$ of the diagonal containing
$\left(p, \bar p\right)\in M\times\bar M$
(here $\bar M$ denotes the manifold
conjugated to $M$).
The diastasis function is the
\K\ potential $D_p$
around $p$ defined by
\[
D_p\left(q\right)=\tilde\Phi \left(q, \bar q\right)+
\tilde\Phi \left(p, \bar p\right)-\tilde\Phi \left(p, \bar q\right)-
\tilde\Phi \left(q, \bar p\right).
\]

If  $d_p:\exp_p\left(V\right)\subset M\rightarrow\R$ denotes  the geodesic distance from $p$
then  one has:
\[
D_p\left(q\right)=d_p\left(q\right)^2+O\left(d_p\left(q\right)^4\right)
\]
and   $D_p=d_p^2$ if and only if $g$ is the flat metric.
We refer the reader to the seminal paper of Calabi  \cite{CalabiIsometric} for more details
and further results on the diastasis function (see also \cite{AndreaDiast}, \cite{AndreaExp}
and \cite{AndreaSymm}).

In \cite{AndreaExp}  it is proven  that there exists an open \ngh\
$S$ of the zero section of the tangent bundle $TM$ of $M$ and a smooth embedding $\nu:S\rightarrow TM$
such that $p\circ \nu =p$, where $p:TM\rightarrow M$ is  the natural projection,  satisfying the following conditions: if one writes
\[
\nu \left(p, v\right)=\left(p, \nu_p\left(v\right)\right),\ \left(p, v\right)\in S
\]
then the diffeomorphism
\[
\nu_p:T_pM\cap S\rightarrow T_pM\cap \nu \left(S\right)
\]
satisfies
\[
\left(d\nu_p\right)_0=\id_{T_pM}
\]
\[
D_p\left( \exp_p \left(\nu_p \left(v\right)\right)\right)=g_p\left(v, v\right),\  \forall v\in T_pM\cap S,
\]
where
$\exp_p:V\subset T_pM\rightarrow M$ denotes the exponential map at $p$
($V$ is a suitable \ngh\ of the origin of  $T_pM$ where the restriction of $\exp_p$
is a diffeomorphism).
 Thus, the smooth map
\[
\E_p:=\exp_p\circ\ \nu_p :T_pM\cap S\rightarrow M
\]
satisfies
\begin{equation}\label{expnup}
\left(d \E_p\right)_0=\id_{T_pM}
\end{equation}
 \begin{equation}\label{expexpdiast}
D_p\left( \E_p \left(v\right)\right)=g_p\left(v, v\right),\  \forall v\in W.
\end{equation}

In analogy
with the exponential at $p$ (which satisfies
$d_p\left( \exp_p \left(v\right)\right)=\sqrt{g_p\left(v, v\right)},$
$\forall v\in V$)
any  smooth map
$\E_p:W\rightarrow M$
from a \ngh\  $W$ of the origin of $T_pM$ into $M$
satisfying (\ref{expnup}) and (\ref{expexpdiast})
will be called a  \emph{diastatic exponential} at $p$.
It is worth pointing out  (see  \cite{boc} for a proof)  that $\exp_p$ is holomorphic if and only if the metric $g$ is flat and it is not hard to see that  the same assertion holds true for  a diastatic
exponential $\E_p$.\\

In this paper we study the  diastatic exponentials for the Hermitian symmetric
spaces of noncompact type (HSSNT) and their compact duals.
The following examples deal with the  rank one case  and it will be our prototypes
for the general case.

\begin{exmp}\rm\label{complexhyp}
Let  $\C H^n= \{z \in \C ^n \, | \, |z |^2=|z_1|^2+\cdots +|z_n|^2<1 \}$
be  the complex hyperbolic space endowed with the hyperbolic metric,
namely this metric $g^{\hyp}$ whose associated \K\ form is given by
$\omega^{\hyp} =-\frac{i}{2}\partial\bar\partial\log \left(1-|z|^2\right)$.
Thus the diastasis function $D_0^{\hyp}:\C H^n \rightarrow \R$ and the exponential map
$\exp_0^{\hyp}:T_0\C H^n \cong \C^n\rightarrow \C H^n$
around the origin $0\in\C^n$
are given respectively  by
\[
D_0^{\hyp}\left(z\right)=- \log \left(1-| z|^2\right)
\]
and
\[
\exp_0^{\hyp} \left(v\right)=\tanh \left(|v|\right)\frac{v}{|v|}, \quad  \exp_0^{\hyp}\left(0\right)=0.
\]
It is then immediate to verify that the map  $\E_0^{\hyp}:T_0\C H^n\rightarrow \C H^n$ given by:
\[
\E_0^{\hyp}\left(v\right)=\sqrt{1- e^{-|v|^2}}\frac{v}{|v|}, \quad \E_0^{\hyp}\left(0\right)=0, \quad v=\left(v_1, \dots v_n\right)
\]
satisfies
$\left(d \E_0^{\hyp}\right)_0=\id_{T_0\C H^n}$ and
\[
D_0^{\hyp}\left( \E_0^{\hyp} \left(v\right)\right)=g^{\hyp}_0\left(v, v\right)=|v|^2, \quad \forall v\in T_0\C H^n =\C^n.
\]
Hence $\E_0^{\hyp}$ is a diastatic exponential at $0$.
Notice that  $ \E_0^{\hyp}$
is  characterized by the fact that
it is direction preserving. More precisely,
if  $F: T_0\C H^n\rightarrow \C H^n$
is a diastatic exponential satisfying
$F\left(v\right)=\lambda \left(v\right)v$,
for some smooth nonnegative   function $\lambda: \C^n\rightarrow \R$,
then $F=\E_0^{\hyp}$.
\end{exmp}

\begin{exmp}\rm\label{complexhyppoly}
 Let $P=\left(\C H^1\right)^{\ell}$  be a polydisk.
If $z_k$,   $k=1, \dots ,\ell$, denotes the complex coordinate in each factor of
  $P$  and  $v=\left(v_1, \dots , v_{\ell}\right)\in T_0P\cong\C^\ell$.
Then the diastasis
 $D_0^{P}:P\rightarrow\R$,
 the exponential map $\exp_0^{P}:T_0P\rightarrow P$
 and a diastatic exponential
 $\E_0^{P}:T_0P\rightarrow P$ at the origin are given respectively by:
\[
D_0^{P}\left(z\right)=- \sum_{k=1}^{\ell}\log \left(1- | z_k |^2\right),
\]
\[
\exp_0^{P} \left(v\right)=\left(\tanh \left(|v_1|\right)\frac{v_1}{|v_1|}, \dots, \tanh \left(|v_{\ell}|\right) \frac {v_\ell} {|v_\ell|}\right), \quad \exp^{\hyp}_0 \left(0\right)=0,
\]
\begin{equation}\label{expdiastpol}
\E_0^{P}\left(v\right)=\left(\sqrt{1- e^{-|v_1|^2}}\frac{v_1}{|v_1|}, \dots, \sqrt{1- e^{-|v_\ell|^2}}\frac{v_\ell}{|v_\ell|}\right),\quad \E_0^{P}\left(0\right)=0.
\end{equation}
\end{exmp}

Let now $M$ be an HSSNT
which we identify with a   bounded symmetric domain of $\C ^n$ centered at the origin $0\in \C^n$
equipped with the hyperbolic   metric $g^{\hyp}$,
namely the \K\ metric whose associated \K\ form (in the irreducible case)
is given by
\[
\omega^{\hyp}= \frac{i}{2g}\de\deb\log K_M.
\]
Here $K_M\left(z, \bar z\right)$  (holomorphic in the first variable
and antiholomorphic in the second one) denotes the reproducing kernel of $M$
and $g$ its genus.
By using the rotational symmetries of $M$
one can show that the diastasis function at the origin
$D^{\hyp}_0:M\rightarrow\R$ is globally defined and  reads as
\[
D^{\hyp}_0\left(z\right)=\frac{1}{g}\log K_M\left(z, \bar z\right),
\]
(see \cite{AndreaDiast} for a proof and  further results on Calabi's function
for HSSNT).
Notice also that, by Hadamard theorem,   the exponential map $\exp^{\hyp}_0:T_0M\rightarrow M$ is a global diffeomorphism.\\

The following theorem which is the first  result of this paper, contains a description of the diastatic exponential for   HSSNT.

\begin{thm}\label{mainteor1}
Let $\left(M, g^{\hyp}\right)$ be an HSSNT.
Then   there exists a globally defined diastatic exponential
$\E^{\hyp}_0:T_0M\rightarrow M$
which is a diffeomorphism and is
uniquely determined by the fact that
 ${\E^{\hyp}_0}_{|T_0P}=\E_0^{P}$
for every polydisk $P\subset M$,  $0\in P$, where $\E_0^{P}$
is given by  (\ref{expdiastpol}).
In particular  ${\E^{\hyp}_0}_{|T_0N}=\E_0^N$ for every complex and totally geodesic submanifold $N\subset M$ through $0$.
\end{thm}

Consider now
the  Hermitian symmetric spaces of compact type (HSSCT).
Let us consider first the compact duals of Examples \ref{complexhyp} and \ref{complexhyppoly}.

\begin{exmp}\rm\label{complexproj}
Let $\C P^n$ be the complex projective space endowed with the Fubini--Study metric $g^{FS}$,
namely the metric whose associated \K\ form is given by
\[
\omega^{FS}=\frac{i}{2}\de \deb \log \left(|Z_0|^2+\dots+|Z_n|^2\right)
\]
for a choice of homogeneous coordinates $Z_0,\dots, Z_n$.
Let $p_0=[1,0,\dots,0]$ and consider the affine chart $U_0=\{Z_0\neq 0\}.$
Thus we have the following inclusions
\begin{equation}\label{inclusions}
\C H^n\subset \C^n\cong U_0\subset \C P^n,\end{equation}
where we are identifying $U_0$ with $\C^n$
via the affine coordinates
\[
U_0\rightarrow\C^n: [Z_0,\dots , Z_N]\mapsto \left(z_1=\frac{Z_1}{Z_0},\dots , z_n=\frac{Z_n}{Z_0}\right).
\]
Under this identification we make no distinction between the point $p_0$ and the origin $0\in \C^n$.
Calabi's  diastasis function $D_0^{FS}:U_0 \rightarrow \R$ around $p_0\equiv 0$ is given by
\[
D_0^{FS}\left(z\right)= \log \left(1 + | z |^2\right).
\]
Observe  that  $D_0^{FS}$ blows up at the points belonging to $\C P^n\setminus U_0$
which is the  cut locus of $p_0$ with respect to the Fubini--Study metric.
We denote this set by $\Cut_0\left(\C P^n\right)$.

It is not hard to verify that the map
\[\E_0^{FS}:T_0\C P^n\rightarrow \C P^n\setminus \Cut_0\left(\C P^n\right)\]
given by
\[
\E_0^{FS}\left(v\right)=\sqrt{e^{|v|^2}-1}\ \frac{v}{|v|}, \quad \E_0^{FS}\left(0\right)=0,
\]
is a diastatic exponential at $0$, namely it satisfies
$\left(d \E_0^{FS}\right)_0=\id_{T_0\C P^n}$ and
\[
D_0^{FS}\left( \E_0^{FS} \left(v\right)\right)=g^{FS}_0\left(v, v\right)=|v|^2, \quad \forall v\in T_0\C P^n .
\]
\end{exmp}
\begin{exmp}\rm\label{complexprojpoly}
Let  $P^*=\left(\C P^1\right)^\ell$
be a (dual) polydisk.
If $z_k$, for $k=1, \dots ,\ell$,  denotes the affine coordinate in each factor of
$P^*$  and  $v=\left(v_1, \dots , v_\ell\right)\in T_0M^*\cong\C^\ell$ then
 it is immediate to see that the diastasis
 $D_0^{P^*}:P^*\rightarrow\R$,
 the exponential map $\exp_0^{P^*}:T_0P^*\rightarrow P^*$
 and a diastatic exponential
 $\E_0^{P^*}:T_0P^*\rightarrow P^*$ at the origin are given respectively by:
\[D_0^{P^*}\left(z\right)= \sum_{k=1}^\ell\log \left(1+ | z_k|^2\right),\]
\[\exp_0^{P^*} \left(v\right)=\left(\tan \left(|v_1|\right)\frac{v_1}{|v_1|}, \dots, \tan \left(|v_\ell|\right) \frac {v_\ell} {|v_\ell|}\right), \quad \exp_0 \left(0\right)=0,
\]
\begin{equation}\label{expdiastpold}
\E_0^{P^*}\left(v\right)=\left(\sqrt{ e^{|v_1|^2}-1}\frac{v_1}{|v_1|}, \dots, \sqrt{e^{|v_\ell|^2}-1}\frac{v_\ell}{|v_\ell|}\right),\quad \E_0^{P^*}\left(0\right)=0.
\end{equation}
\end{exmp}

Given an arbitrary   HSSNT $M$ of genus $g$ let denote by $M^*$ its compact dual  equipped with  the Fubini--Study metric  $g^{FS}$, namely the pull-back of the Fubini--Study metric of $\C P^N$
via the Borel--Weil embedding $M^*\rightarrow\C P^N$  (see \cite{sympdual} for details).
Let $0\in M^*$ be a fixed point and denote by $\Cut_0\left(M^*\right)$ the cut locus of $0$ with respect to
the Fubini--Study metric.
In the irreducible case the \K\ form $\omega^{FS}$ associated to $g^{FS}$
is given (in the {\em affine} chart  $M^*\setminus\Cut_0\left(M^*\right)$) by
\[
\omega^{FS}=\frac{i}{2g}\de\deb\log K_{M^*},
\]
where
\begin{equation}\label{KMD}
K_{M^*}\left(z, \bar z\right)=1/K_M\left(z, - \bar z\right).
\end{equation}
We call $K_{M^*}$ the {\em reproducing kernel} of $M^*$.
Notice that $K_{M^*}$  is the weighted Bergman  kernel for the (finite dimensional) complex Hilbert space consisting of  holomorphic functions $f$ on
$M^*\setminus\Cut_0\left(M\right)\subset M^*$ such that $\int_{M^*\setminus\Cut_0\left(M\right)}|f|^2\left(\omega^{FS}\right)^n < \infty$ (see \cite{AndreaDiast} and also \cite{me} for a nice characterization of symmetric spaces in terms of $K_{M^*}$).
Notice that when  $M=\C H^n$ then $g=n+1$, $K_M\left(z, \bar z\right)=\left(1-|z|^2\right)^{-\left(n+1\right)}$,
$K_{M^*}\left(z, \bar z\right)=\left(1+|z|^2\right)^{n+1}$  and the Borel--Weil embedding is the identity of $\C P^n$.\\

Observe that,  as in the previous examples,  $D_0^{FS}$ is globally defined in
$M^*\setminus \Cut_0\left(M^*\right)$ (see \cite{tas} for a proof)
and it blows up at the points in $\Cut_0\left(M^*\right)$.
Moreover
\[
D_0^{FS}\left(z\right)=\frac{1}{g}\log K_{M^*}\left(z, \bar z\right), \quad z\in M^*\setminus\Cut_0\left(M^*\right).
\]

Furthermore (see e.g. \cite{wolf})
$M^*\setminus \Cut_0\left(M^*\right)$ is globally  biholomorphic to
$T_0M$ and if $0$ denote the origin of $M$   one has the following inclusions (analogous of (\ref{inclusions}))
\begin{equation}\label{inclusionsgen}
M\subset T_0M=T_0M^*\cong M^*\setminus \Cut_0\left(M^*\right)\subset M^* .
\end{equation}

We are now in the position to state our second result which is the dual counterpart of Theorem \ref{mainteor1}.

\begin{thm}\label{mainteor1d}
Let $\left(M^*, g^{FS}\right)$
 be an HSSCT.
Then  there exists a globally defined diastatic exponential
$\E_0^{FS}:T_0M^*\rightarrow M^*\setminus \Cut_0\left(M^*\right)$
which is
uniquely determined by the fact that for
every (dual) polydisk $P^*=\left(\C P^1\right)^s\subset M^*$
its restriction to $T_0P^*$
equals the map   $\E_0^{P^*}$ given by
(\ref{expdiastpold}).
In particular ${\E_0^{FS}}_{|T_0N^*}=\E_0^{N^*}$
for every complex and totally geodesic submanifold $N^*\subset M^*$ through $0$.
\end{thm}

The key ingredient for the proof of Theorem \ref{mainteor1} and Theorem \ref{mainteor1d} is the theory of Hermitian positive Jordan
triple systems (HPJTS).
In \cite{sympdual} this theory has been the main tool  to study the link between the  symplectic geometry of an  Hermitian symmetric space $\left(M, \omega^{\hyp}\right)$ and its dual $\left(M^*, \omega^{FS}\right)$ where $\omega^{\hyp}$ (resp. $\omega^{FS}$) is the \K\ form associated to  $g^{\hyp}$  (resp.  $g^{FS}$).
The main result proved there, is the following theorem.

\begin{thm}\label{thmsympdual}
Let $M$ be an HSSNT and $B\left(z, w\right)$ its associated Bergman operator (see next section).
Then the map
\begin{equation}\label{sympdualmap}
\Psi_M :M\rightarrow  M^*\setminus \Cut_0\left(M^*\right),\quad z\mapsto B\left(z, z\right)^{-\frac{1}{4}}z
\end{equation}
is a global real analytic diffeomorphism such that
\[
\Psi_M^*\omega_0=\omega^{\hyp}
\]
\[
\Psi_M^*\omega^{FS}=\omega_0,
\]
where $\omega_0$ is the flat \K\ form on $T_0M$.
Moreover,
for every complex and totally geodesic submanifold $N\subset M$ one has
${\Psi_M}_{|N}=\Psi_N$.
\end{thm}

Here  $\omega_0$ denotes the  \K\ form
 on $M$ obtained by the restriction of the flat \K\ form on $T_0M=\C^n$.
The map  $\Psi_M$ was christened in \cite{sympdual} as the  \emph {symplectic duality}.
The unicity of this map and an alternative  proof of Theorem \ref{thmsympdual}
can be found in \cite{unicdual}.

The following theorem which represents our  third result  provides a geometric interpretation of the symplectic duality  map in terms of diastatic exponentials.

\begin{thm}\label{mainteor3}
Let $M$ be a HSSNT and $M^*$ be  its compact dual.
Then the symplectic duality map can be written
as
\[
\Psi_M =\E_0^{FS}\circ \left({\E_0^{\hyp}}\right)^{-1}:M\rightarrow M^*\setminus \Cut_0 \left(M^*\right),
\]
where $\E_0^{\hyp}:T_0M\rightarrow M$ and
$\E_0^{FS}:T_0M^*\rightarrow M^*\setminus \Cut_0\left(M^*\right)$
are  the diastatic exponentials at $0$ of  $M$ and $M^*$ respectively.
\end{thm}

Our fourth result is the following theorem
which shows  that
the \lq\lq algebraic manipulation''  (\ref{KMD}) which allows us to pass from
$K_M$ to $K_{M^*}$  can be realized via the symplectic duality map.

\begin{thm}\label{mainteor4}
Let $K_M$ be the reproducing kernel for an HSSNT and let $K_M^*$
be its dual.
Then
\[
K_{M^*}\circ\Psi_M =K_M,
\]
where $\Psi_M :M\rightarrow M^*\setminus \Cut_0\left(M^*\right)$
is the symplectic duality map.
\end{thm}

\vskip 0.3cm
The paper contains another section,
where,  after recalling
some standard facts about HSSNT and HPJTS,
we prove Theorem \ref{mainteor1}, Theorem \ref{mainteor1d}, Theorem \ref{mainteor3}
and
Theorem \ref{mainteor4}.

\section{hpjts and the proofs of the  main results}
We refer the reader  to \cite{roos} (see also \cite{loos}) for more details
of  the material on Hermitian positive Jordan triple systems.

\subsection{Definitions and notations}
An Hermitian Jordan triple system is a  pair $\left({\mathcal M},
\{ ,  ,\}\right)$, where ${\mathcal M}$ is a complex vector space and $\{
,  ,\}$ is a map
\[
\{ ,  ,\}:{\mathcal M}\times {\mathcal M}\times {\mathcal M} \rightarrow {\mathcal M}
\]
\[
\left(u, v, w\right)\mapsto \{u, v, w\}
\]
which is ${\C}$-bilinear and symmetric in $u$ and $w$, ${\C}$-antilinear in $v$ and such that the following {\emph Jordan identity} holds:
\[
\{x, y, \{u, v, w\}\}-\{u, v, \{x, y, w\}\}= \{\{x, y, u\}, v,
w\}-\{u, \{v, x, y\}, w\}.
\]
For $x,y,z \in \M$ considered  the following operator
\[
T\left(x,y\right)z =\left\{  x,y,z\right\} 
\]
\[
Q\left(x,z\right)y =\left\{  x,y,z\right\}  
\]
\[
Q\left(x,x\right) =2Q\left(x\right)\label{D3}\\
\]
\[
B\left(x,y\right) =\operatorname{id}_{\mathcal M}-T\left(x,y\right)+Q\left(x\right)Q\left(y\right). \label{D4}%
\]
The operators $B\left(x,y\right)$ and $T\left(x,y\right)$ are $\mathbb{C}$-linear, the operator
$Q\left(x\right)$ is $\mathbb{C}$-antilinear. $B\left(x,y\right)$ is called the \emph {Bergman operator}.
For $z\in V$, the \emph{odd powers} $z^{\left(2p+1\right)}$ of $z$ in the Jordan triple
system $V$ are defined by
\[
z^{\left(1\right)}=z \qquad z^{\left(2p+1\right)}=Q\left(z\right)z^{\left(2p-1\right)}. \label{D6}%
\]
An Hermitian Jordan triple system is called  {\emph positive} if the Hermitian form
\[
\left(  u\mid v\right)  =\tr T\left(u,v\right) \label{D5}%
\]
is positive definite. An element $c \in \M$ is called \emph {tripotent} if
$\{c,c,c\}=2c$. Two tripotents $c_1$ and $c_2$ are called \emph {(strongly)
orthogonal} if $T\left(c_1, c_2\right)=0$.

\subsection{HSSNT associated to HPJTS}
M. Koecher (\cite{Koecher1}, \cite{Koecher2}) discovered that to every HPJTS
$\left(\M, \{ ,  ,\}\right)$ one can associate an Hermitian symmetric
space of noncompact type, i.e. a bounded symmetric domain $M$
centered at the origin $0\in \M$. The domain $M$ is defined as the connected component containing the origin of   the set of all $u\in {\M}$ such that $B\left(u, u\right)$ is positive definite with respect to the Hermitian form $\left(u, v\right)\mapsto \tr T\left(u, v\right)$. \emph{We will always consider such a domain in its (unique up to linear isomorphism) circled realization.}
The reproducing kernel $K_M$ of $M$ is given by
\begin{equation}\label{KOB}
K_M\left(z, \bar z\right)= \det B\left(z,z\right)
\end{equation}
and so when $M$ is irreducible
\[
\omega^{\hyp}= -\frac{i}{2g} \partial \bar \partial \log\det B.
\]

The HPJTS $\left({\M}, \{ ,  ,\}\right)$ can be recovered by its
associated HSSNT $M$ by defining ${\M}=T_0 M$ (the tangent space to the origin of $M$) and
\begin{equation}\label{trcurv}
\{u, v, w\}=-\frac{1}{2}\left(R_0\left(u, v\right)w+J_0R_0\left(u, J_0v\right)w\right),
\end{equation}
where $R_0$ (resp. $J_0$) is the curvature tensor of the Bergman metric (resp. the complex structure) of $M$ evaluated at the origin.
The reader is referred   to  Proposition III.2.7 in  \cite{Bertram} for the proof of (\ref{trcurv}). For more informations on the correspondence between
$HPJTS$ and $HSSNT$ we  refer also  to p. 85 in Satake's book \cite{satake}.

\subsection{Totally geodesic submanifolds of HSSNT}
In the proof of our theorems  we  need the following result.
\begin{prop} \label{sub}
Let $M$ be a HSSNT and let ${\M}$ be its associated $HPJTS$. Then there exists a one to one
correspondence between (complete) complex totally geodesic
submanifolds through the origin and sub-$HPJTS$ of ${\M}$. This correspondence
sends $T \subset M$ to  ${\mathcal T} \subset {\M}$, where
${\mathcal T}$ denotes the HPJTS associated to $T$.
\end{prop}

\subsection{Spectral decomposition and Functional calculus}
Let $\M$ be a HPJTS. Each element $z\in \M$ has a unique \emph{spectral decomposition}%
\[
z=\lambda_{1}c_{1}+\cdots+\lambda_{s}c_{s}\qquad\left(0<\lambda_{1}<\cdots
<\lambda_{s}\right), \label{D7}%
\]
where $\left(c_{1},\ldots,c_{s}\right)$ is a sequence of pairwise
orthogonal tripotents and the $\lambda_j$ are real number called eigenvalues of $z$. For every $z \in \M$ let $\max\{z\}$ denote the largest eigenvalue of $z$, then $\max\{\cdot \}$ is a norm on $\M$ called the \emph{spectral norm}. The HSSNT $M$ associated to $\M$
is the open unit ball in $\mathcal M$ centered at the origin (with respect the spectral norm $M$),
i.e.,
\begin{equation}\label{Mball}
M=\{z=\sum_{j=1}^s\lambda_jc_j \ |\ \max\{z\}= \max_j\{ \lmb _j\}<1\}
\end{equation}

Using the spectral decomposition, it is possible to associate to an \emph{odd} function
$f:\mathbb{R}\rightarrow \mathbb{C}$ a map  $F:\M \rightarrow \M$
as follows. Let $z\in \M$
and let
\[
z=\lambda_{1}c_1+\cdots+\lambda_{s}c_s,\quad 0<\lambda_{1}<\cdots<\lambda
_{s}
\]
be the spectral decomposition of $z$. Define the map $F$
by
\begin{equation}\label{associatedfunction}
F\left(z\right)=f\left(\lambda_1\right)c_1+\cdots+f\left(\lambda_s\right)c_s.
\end{equation}
If $f$ is continuous, then $F$ is continuous. If
\[
f\left(t\right)=\sum_{k=0}^{N}a_{k}t^{2k+1}%
\]
is a polynomial, then $F$ is the map defined by%
\[
F\left(z\right)=\sum_{k=0}^{N}a_{k}z^{\left(2k+1\right)}\qquad\left(z\in \M\right). \label{F02}%
\]
If $f$ is analytic, then $F$ is real-analytic. If $f$ is given near $0$ by%
\[
f\left(t\right)=\sum_{k=0}^{\infty}a_{k}t^{2k+1},
\]
then $F$ has the Taylor expansion near $0\in V$:%
\[
F\left(z\right)=\sum_{k=0}^{\infty}a_{k}z^{\left(2k+1\right)}. \label{F03}%
\]

\begin{exmp}\label{polysp}\rm
Let $P=\left(\C H^1\right)^\ell \subset \left(\C^\ell, \{,,\}\right)$ be the polydisk embedded in is its associated HPJTS $\left(\C^\ell,\{,,\}\right)$. Define $\tilde c_j= \left(0,\dots,0, e^{i \theta_j},0,\dots,0\right),$ $1\leq j \leq \ell$. The $\tilde{c}_j$ are mutually strongly orthogonal tripotents. Given $z=\left(\rho_1 e^{i \theta_1}, \dots, \rho_\ell e^{i \theta_\ell}\right) \in \left(\C H^1\right)^\ell,$ $z\neq 0,$ then up to a permutation of the coordinates, we can assume $0 \leq \rho_1 \leq \rho_2 \leq \dots \leq \rho_\ell$. Let $i_1,$ $1 \leq i_1 \leq \ell,$ the first index such that $\rho_{i_1} \neq 0$ then we can write
$$
z=\rho_{i_1}\left(\tilde{c}_{i_1} + \dots + \tilde{c}_{i_2-1}\right) + \rho_{i_2}\left(\tilde{c}_{i_2} + \dots + \tilde{c}_{i_3-1}\right) + \dots + \rho_{i_s}\left(\tilde{c}_{i_s} + \dots + \tilde{c}_{i_{s+1}-1}\right)
$$
with $0 < \rho_{i_1} < \rho_{i_2} < \dots < \rho_{i_s}=\rho_\ell$ and $i_{s+1}=\ell+1$. The $c_j$'s, defined by $c_j =\tilde{c}_{i_j}+ \dots +\tilde{c}_{i_{j+1}-1}$, are still mutually strongly orthogonal tripotents and $z=\lmb_1 c_1 + \dots + \lmb_s c_s$ with $\lmb_j=\rho_{i_j},$ is the spectral decomposition of $z$. So the diastatic exponential given in (\ref{expdiastpol}) can be written as
\[
\E_0^{P}\left(z\right)=\left(\sqrt{1- e^{-|z_1|^2}}\frac{z_1}{|z_1|}, \dots, \sqrt{1- e^{-|z_\ell|^2}}\frac{z_\ell}{|z_\ell|}\right)= \sum^s_{j=1} \left(1-e^{-\lambda_j^2}\right) ^\frac{1}{2} c_j
\]
and $\E_0^{P}\left(0\right)=0$.
\end{exmp}

We are now in the position to prove our main results.
In all the following proofs we can assume,
without loss of generality, that $M$ is  irreducible.  Indeed, in  the reducible case the Bergman operator
is the product of the Bergman operator  of each factor and therefore the same holds  true for the diastatic exponential and for the symplectic duality map.

\subsection{Proof of Theorem \ref{mainteor1}}
Consider the odd smooth function $f:\R\rightarrow\R$
defined by
\[f\left(t\right)=\left(1-e^{-t^2}\right)^\frac{1}{2}\frac{t}{|t|}, \ f\left(0\right)=0\]
and the map
$F:T_0M \rightarrow M \subset T_0M$
associated to $f$ by (\ref{associatedfunction}),
namely
\begin{equation}\label{expdiastgen}
F\left(z\right)=\sum_{j=1}^s \left(1-e^{-\lmb^2_j}\right)^{\frac{1}{2}}c_j,
\end{equation}
where  $z=\lambda_1 c_1+ \dots +\lambda_s c_s$ is the spectral decomposition of $z \in M$. Notice that $F\left(T_0M\right) \subset M$ by (\ref{Mball}).
We will show that   $\E_0^{\hyp}:=F$ is indeed a diastatic exponential at the origin
for $M$ satisfying the conditions of Theorem \ref{mainteor1}. It is easy to see that $\E_0^{\hyp}$ is injective and $\left(d\E_0^{\hyp}\right)_0=id_{T_0M}.$
Thus,  it remains to show that
$D_0^{\hyp}\left(\E_0^{\hyp} \left(z\right)\right)=  g_0^{\hyp}\left(z,z\right)$.
In order to prove this  equality observe that  (see \cite{roos} for a proof)
\begin{equation}\label{BZZ}
B\left(z,z\right)c_j=\left(1-\lambda_j^2\right)^2c_j, \quad j=1,\dots ,s,
\end{equation}
\[
\det B\left(z,z\right) = \prod_{j=1}^s\left(1-\lmb_j^2\right)^g,
\]
\[
g^{\hyp}_0\left(z, z\right)= \frac{1}{g} \tr T\left(z,z\right)= \sum_{j=1}^s \lmb_j^2.
\]
Thus (\ref{KOB}) yields,
\begin{equation}\label{doz}
D^{\hyp}_0\left(z\right)=-\frac{1}{g} \log \det B\left(z,z\right)=- \log \prod_{j=1}^s\left(1-\lmb_j^2\right)
\end{equation}
and so
\[
D^{\hyp}_0\left(\E_0^{\hyp} \left(z\right)\right)= - \log \prod_{j=1}^s \left[1-\left(1-e^{-\lmb_j^2}\right)\right] = \sum_{j=1}^s\lmb_j^2 = g_0^{\hyp}\left(z,z\right),
\]
namely the desired equality.
In order to prove the second part of the theorem let $P \subset M$ be a polydisk through the origin.
Thus equality  ${\E^{\hyp}_0}_{|T_0P}=\E_0^{P}$ follows by Proposition \ref{sub}, Example \ref{polysp}
and formula (\ref{expdiastgen}). Moreover $\E^{\hyp}_0$ is determined by its restriction to polydisks since it is well-known that  $\forall z \in T_0M$ there exists a polydisk $P \subset M$ such that $0 \in P$ and $z \in T_0P$ (see, e.g.
\cite{helgason} and also
\cite{Polydisc}).

\subsection{Proof of Theorem \ref{mainteor1d}}
Let $z=\lambda_1 c_1 + \dots + \lambda_s c_s$ be a spectral decomposition of $z \in M^* \setminus \Cut_0\left(M^*\right) \cong T_0M$. In analogy with the compact case  one has
\[
B\left(z,-z\right)c_j=\left(1+\lambda_j^2\right)^2c_j
\]
\[
\det B\left(z,-z\right) = \prod_{j=1}^s\left(1+\lmb_j^2\right)^{g}.
\]
\[
g^{FS}_0\left(z, z\right)= \lmb_j^2.
\]
Thus,  by (\ref{KMD}),  Calabi's diastasis function at the origin for $g^{FS}$
is given by:
\begin{equation}\label{dozd} \begin{split}
D_0^{FS}\left(z\right)&= - \frac{1}{g} \log K_{M^*}\left(z, \bar z\right) = \frac{1}{g} \log [K_M\left(z, -\bar z\right)] = \frac{1}{g} \log[\det B\left(z,-z\right)]  \\
& = \frac{1}{g} \log \prod_{j=1}^s\left(1+\lmb_j^2\right)
\end{split}\end{equation}
Define
$
\E_0^{FS}:T_0M^* \cong T_0M \f  M^* \setminus \Cut_0\left(M^*\right) \cong T_0M
$
as the map associated to the real function
$
f^*\left(t\right)=\left(e^{t^2}-1\right)^\frac{1}{2} \frac{t}{|t|}
$ by (\ref{associatedfunction}),
namely
\begin{equation}\label{expdiastgend}
\E_0^{FS}\left(z\right) = \sum_{j=1}^s \left(e^{\lmb^2_j}-1\right)^{\frac{1}{2}}c_j.
\end{equation}
Thus, following the same line of the proof of Theorem \ref{mainteor1},
one can show that $\E_0^{FS}$ is  the diastatic exponential at $0$ uniquely determined by its restriction to polydisks.

\subsection{Proof of Theorem \ref{mainteor3}}
By (\ref{sympdualmap}) and (\ref{BZZ})
\begin{equation}\label{PSPEC}
\Psi_M\left(z\right)=B\left(z,\bar z\right)^{-\frac{1}{4}}\left(z\right)= \frac{\lmb_j}{\left(1-\lmb_j^2\right)^{\frac{1}{2}}}c_j
\end{equation}
By the very definition of the diastatic exponential $\E_0^{\hyp}$ for the hyperbolic metric
its inverse
$\left({\E_0^{\hyp}}\right)^{-1}:M\rightarrow {T_0M}$
read as:
\[\left({\E_0^{\hyp}}\right)^{-1}\left(z\right)=\sum_{j=1}^s \left(-\log\left(1-\lmb^2_j\right)\right)^{\frac{1}{2}}c_j,\]
Then,  by (\ref{expdiastgend}) and (\ref{PSPEC}),
\[
\E_0^{FS} \circ \left({\E_0^{\hyp}}\right)^{-1}\left(z\right)=\Psi_M\left(z\right)
\]
and this concludes the proof of Theorem \ref{mainteor3}.

\subsection{Proof of Theorem \ref{mainteor4}}
Since $D_0^{\hyp}=\frac{1}{g}\log K_M$ and
$D_0^{FS}=\frac{1}{g} \log K_{M^*},$
equation $K_{M^*}\circ\Psi_M =K_M$
is equivalent to
$D_0^{FS}\circ\Psi_M =D^{\hyp}_0$ which is a straightforward consequence of (\ref{doz}), (\ref{dozd}) and (\ref{PSPEC}).

\end{document}